\documentclass[12pt]{article}
\usepackage{setspace}
\usepackage{amsmath,amssymb,amsthm}
\usepackage{mathrsfs}
\usepackage{geometry}
\usepackage{enumitem}
\usepackage{environ}
\usepackage{setspace}
\usepackage{titlesec}
\usepackage[all,cmtip]{xy}
\usepackage{tikz}
\usetikzlibrary{matrix,arrows}

\NewEnviron{prf}[1][]{\begin{proof}[\bf #1Proof]\BODY\end{proof}}{}
\NewEnviron{slt}[1][]{\begin{proof}[\bf #1	解]\BODY\end{proof}}{}
\newtheorem{definition}{Definition}[subsection]
\newtheorem{theorem}[definition]{Theorem} 
\newtheorem{lemma}[definition]{Lemma}
\newtheorem{proposition}[definition]{Proposition}
\newtheorem{remark}[definition]{Remark}
\newtheorem{example}[definition]{Example}
\newtheorem{corollary}[definition]{Corollary}
\newtheorem{conjecture}[definition]{Conjecture}

\setlist[enumerate,1]{label=(\roman*)}
\setlist[enumerate,2]{label=(\alph*)}

\title{\textbf{Dynamical Mordell--Lang conjecture for totally inseparable liftings of Frobenius}}
\author{She YANG}
\date{}

\begin{document}
\begin{spacing}{1.25}

\maketitle

\begin{abstract}
We prove that if $K$ is a complete algebraically closed non-archimedian valuation field of positive characteristic and $f$ is an endomorphism of $\mathbb{P}_{K}^{N}$ which is totally inseparable and behaves as the Frobenius on the special fiber, then $f$ satisfies the dynamical Mordell--Lang (DML) property. We also discuss some corollaries and generalizations.
\end{abstract}

\section{Introduction}
In this paper, as a matter of convention, every variety is assumed to be integral but the closed subvarieties can be reducible. We denote $\mathbb{N}=\mathbb{Z}_{+}\cup\{0\}$. An arithmetic progression is a set of the form $\{mk+l|\ k\in\mathbb{N}\}$ for some $m,l\in\mathbb{N}$. For a set $X$ endowed with a self-map $f$ and for $x\in X$, we define the orbit of $x$ under $f$ as $\mathcal{O}_{f}(x):=\{f^{n}(x)|\ n\in\mathbb{N}\}$. A subset $Y\subseteq X$ is called $f$-periodic if $f^{n}(Y)\subseteq Y$ for some $n\in\mathbb{Z}_{+}$.

Let $K$ be a field. Let $X$ be a quasi-projective variety over $K$ and let $f:X\rightarrow X$ be a $K$-endomorphism. We say that $f$ satisfies the DML property if for every $x\in X(K)$ and every closed subvariety (i.e. reduced closed subscheme) $V$ of $X$, the set $\{n\in\mathbb{N}|\ f^{n}(x)\in V(K)\}$ is a finite union of arithmetic progressions. The dynamical Mordell--Lang conjecture in characteristic 0, which is one of the core problems in the field of arithmetic dynamics, is stated as follows. It was first proposed in $\cite{GT09}$ under the influence of S.-W. Zhang (see $\cite[\mathrm{top}\ \mathrm{of}\ \mathrm{p.}\ 306]{GT08}$).

\begin{conjecture}$\mathrm{(Ghioca}$--$\mathrm{Tucker)}$
If $\mathrm{char}K=0$, then every endomorphism $f$ of every quasi-projective variety satisfies the DML property.
\end{conjecture}

Many works toward this conjecture have been done. For example, there are $\cite{BGHKST13}$, $\cite{BGKT12}$, $\cite{Fak14}$, $\cite{GT09}$, $\cite{GTZ08}$, $\cite{GTZ12}$ and $\cite{GX20}$. Two notable cases are the followings:
\begin{enumerate}
\item
If $f$ is an \'etale endomorphism, then $f$ satisfies the DML property. See $\cite[\mathrm{Theorem}\ 1.3]{BGT10}$.
\item
If $X=\mathbb{A}^{2}$ and $f$ is an endomorphism of $X$, then $f$ satisfies the DML property. See $\cite{Xie17}$.
\end{enumerate}

However, the naive analogue of Conjecture 1.0.1 is false when char$K>0$. For an example, see $\cite[\mathrm{Example}\ 3.4.5.1]{BGT16}$. As a result, it is natural to ask two questions:
\begin{enumerate}
\item
What is the form of the set $\{n\in\mathbb{N}|\ f^{n}(x)\in V(K)\}$ when char$K>0$?
\item
In which case does the endomorphism $f$ satisfy the DML property?
\end{enumerate}

For the first question, Ghioca and Scanlon proposed a conjecture $\cite[\mathrm{Conjecture}\ 13.2.0.1]{BGT16}$, which is known as the dynamical Mordell--Lang conjecture in positive characteristic. One can refer to $\cite{CGSZ21}$ for the latest progress on this conjecture.

On the other hand, for the second question, $\cite{Xie}$ (and also Ghioca and Scanlon) guessed that most dynamic systems should satisfy the DML property and the counterexamples often involve some group actions. This paper is an attempt towards this opinion.

To our knowledge, when char$K>0$, there are only two nontrivial known cases in which $f$ satisfies the DML property:
\begin{enumerate}
\item
If $X$ is a semiabelian variety defined over a finite field and $\Phi$ is an algebraic group endomorphism satisfies $f(\Phi)=0$ in which $f\in\mathbb{Z}[x]$ is a polynomial with leading coefficient 1 such that the nonzero roots of $f$ in $\mathbb{C}$ are distinct, then $\Phi$ satisfies the DML property. See $\cite[\mathrm{Proposition}\ 13.3.0.2]{BGT16}$ and $\cite[\mathrm{middle}\ \mathrm{of}\ \mathrm{p.}\ 671]{CGSZ21}$.
\item
If $X$ is a projective surface and $f$ is an automorphism satisfying $\lambda_{1}(f)>1$ or $X=\mathbb{A}^{2}$ and $f$ is a birational endomorphism satisfying $\lambda_{1}(f)>1$, then $f$ satisfies the DML property. The notion $\lambda_{1}(f)$ stands for the first dynamical degree of $f$. See $\cite[\mathrm{Theorem}\ 1.4]{Xie}$ and $\cite{Xie14}$.
\end{enumerate}

The main theorem of this paper is the following Theorem 1.0.2. We will also prove a more general statement at the end of this paper. See Proposition 5.2.6 and Remark 5.2.7. They are about the lifting of Frobenius on the projective space. We mention that there are several works in $\cite{Hru01}$, $\cite[\mathrm{Subsection}\ 7.3]{MS14}$, $\cite{PR04}$ and $\cite{Xie18}$ about the lifting of Frobenius towards the dynamical Manin--Mumford conjecture and the dynamical Mordell--Lang conjecture for coherent backward orbits. Our result is towards the dynamical Mordell--Lang conjecture for forward orbits.

\begin{theorem}
Let $K$ be a complete algebraically closed non-archimedian valuation field of characteristic $p>0$. Then the endomorphism $f:\mathbb{P}_{K}^{N}\rightarrow\mathbb{P}_{K}^{N}$,
$$\begin{bmatrix} x_{0} \\ x_{1} \\ \vdots \\ x_{N} \end{bmatrix} \mapsto \begin{bmatrix} \sum\limits_{i=0}^{N}a_{0i}x_{i}^{q}+g_{0}(x_{0}^{p},\cdots,x_{N}^{p}) \\ \sum\limits_{i=0}^{N}a_{1i}x_{i}^{q}+g_{1}(x_{0}^{p},\cdots,x_{N}^{p}) \\ \vdots \\ \sum\limits_{i=0}^{N}a_{Ni}x_{i}^{q}+g_{N}(x_{0}^{p},\cdots,x_{N}^{p}) \end{bmatrix}$$
satisfies the DML property, in which $A=(a_{ij})_{(N+1)\times(N+1)}\in GL_{N+1}(\mathcal{O}_{K})$, $q$ is a power of $p$ and $g_{0},\cdots,g_{N}\in\mathfrak{m}_{K}[x_{0},\cdots,x_{N}]$ are homogeneous polynomials of degree $=\frac{q}{p}$.
\end{theorem}

In fact, Theorem 1.0.2 is a generalized version of Theorem 1.0.3 below.

\begin{theorem}
Let $K$ be a local field of positive characteristic. Let $k$ be the residue field of $\mathcal{O}_{K}$ which is finite. Let $f_{\mathcal{O}_{K}}:\mathbb{P}_{\mathcal{O}_{K}}^{N}\rightarrow\mathbb{P}_{\mathcal{O}_{K}}^{N}$ be an $\mathcal{O}_{K}$-morphism which satisfies:\\
1. $f_{\mathcal{O}_{K}}^{*}(\Omega_{\mathbb{P}_{\mathcal{O}_{K}}^{N}/\mathcal{O}_{K}})\rightarrow\Omega_{\mathbb{P}_{\mathcal{O}_{K}}^{N}/\mathcal{O}_{K}}$ is the zero map.\\
2. $f_{0}(=f_{\mathcal{O}_{K}}\times_{\mathcal{O}_{K}}k):\mathbb{P}_{k}^{N}\rightarrow\mathbb{P}_{k}^{N}$ is some $\mathrm{Frob}_{q}$ ($q$ is a power of the prime $\mathrm{char}k$).\\
Then the $\overline{K}$-morphism $f_{\overline{K}}:\mathbb{P}_{\overline{K}}^{N}\rightarrow\mathbb{P}_{\overline{K}}^{N}$ induced by $f_{\mathcal{O}_{K}}$ satisfies the DML property.
\end{theorem}

\begin{remark}
The two conditions in Theorem 1.0.3 could be weakened as $f_{\mathcal{O}_{K}}^{m}$ and $f_{0}^{m}$ satisfying the corresponding properties for some positive integer $m$ because $f_{\overline{K}}^{m}$ satisfies the DML property implies that $f_{\overline{K}}$ also satisfies.
\end{remark}

Since the proof of Theorem 1.0.3 is the technical heart of this paper, we will discuss our strategy towards it in more detail. Firstly, we reduce Theorem 1.0.3 to Theorem 2.1.2 and assume the local field $K=k((t))$ in there without loss of generality (see Remark 2.1.4). Next, we base change the data from $K=k((t))$ to $L=\bar{k}((t))$ in order to make use of our fundamental tool, that is, the jet schemes. 

Now, we have an integral closed subvariety $X_{L}\subseteq\mathbb{P}_{L}^{N}$ and a point $x\in\mathbb{P}_{L}^{N}(L)$ such that $\mathcal{O}_{f_{L}}(x)\cap X_{L}$ is dense in $X_{L}$. We want to show that $X_{L}$ is $f_{L}$-periodic. We construct a model $\mathcal{X}$ of $X_{L}$ over $S=\mathrm{Spec}(\bar{k}[[t]])$. We note that $f_{S}=f_{\mathcal{O}_{K}}\times_{\mathcal{O}_{K}}S$ is attracting at each $\bar{k}$-point $x_{0}$ in the special fiber of $\mathbb{P}_{S}^{N}$. Inspired by hyperbolic dynamical system, we construct a lifting $x\in\mathbb{P}_{L}^{N}(L)=\mathbb{P}_{S}^{N}(S)$ of $x_{0}$ as an analogy of the ``unstable manifold" through $x_{0}$. Under our assumption, if $x_{0}$ lies in the special fiber of $\mathcal{X}$, we may ask $x\in\mathcal{X}(S)$. Our construction is based on the construction of jet schemes.

The jet schemes is an algebraic-geometrically analogue of the jet spaces in complex geometry. It was firstly introduced in $\cite{Ros13}$ in order to give an algebraic-geometrically proof of the classical Mordell--Lang conjecture in positive characteristic. Just the same as in complex geometry, a point in the $n$th jet scheme represents an ``$n$th infinitesimal direction" on the initial scheme. For a rigorous statement, see Proposition 2.3.2. The author would like to regard this property as the core property of jet schemes. The other properties mentioned in Subsection 2.3 are considered as the evidence to show that jet schemes are indeed well-behaved as \emph{schemes} although they are constructed by Weil restriction.

For the critical schemes introduced in Section 3, there is a vivid explanation. We regard a point in the $n$th critical scheme $\mathrm{Crit}^{n}(\mathcal{X},\mathbb{P}_{S}^{N})$ as an ``unstable $n$th direction based on a point in $\mathcal{X}$". The word ``unstable" means that it can be obtained by an $n$th push-forward of the endomorphism $f_{S}$. Hence the image $\mathrm{Exc}^{n}(\mathcal{X},\mathbb{P}_{S}^{N})$ stands for the set of the points in $\mathcal{X}$ which admits an unstable $n$th direction and the abstract statement $\mathrm{Exc}^{n}(\mathcal{X},\mathbb{P}_{S}^{N})=\mathcal{X}$ in Proposition 3.2.1 can be read as each point in $\mathcal{X}$ admits a such direction.

By the way, we would like to mention that the ambient variety in Theorem 1.0.3 can be extended from the projective space to arbitary smooth projective variety using $\cite[\mathrm{Theorem}\ 6.1]{Xie18}$. However, we deliberately focus on the case of projective space. This is not only because we can highlight the principal line of our proof in this way, but the author is more fond of elementary statements and would like to consider Theorem 1.0.2 as the main result of this paper.

~

At the end of the Introduction, we describe the structure of this paper. We will do a reduction (reduce Theorem 1.0.3 to Theorem 2.1.2) and introduce the jet schemes, which is the fundamental tool of our method in Section 2. After that, we will define the critical schemes and use them to deduce a lifting proposition in Section 3. We will finish the proof of Theorem 1.0.3 in Section 4. Lastly, we will propose a corollary and prove Theorem 1.0.2 and its generalization in Section 5.

~

\textbf{Acknowledgement.} I am very grateful to my advisor Junyi Xie who suggested me using the technic of jet schemes in $\cite{Ros13}$ to study the dynamical Mordell--Lang conjecture in positive characteristic. I am also indebted to him for lots of helpful conversations during the preparation of this paper. I am very grateful to Jason P. Bell, Dragos Ghioca and Thomas J. Tucker who examined the earlier version of this article and put forward many helpful suggestions. Moreover, I would like to thank my classmate Xiangqian Yang for some useful discussions.

The author is supported by NSFC Grant (No. 12271007).

\section{Preparations}
We will reduce Theorem 1.0.3 to a geometric version in Subsection 2.1. Then we will introduce the definition and properties of jet schemes. The references of jet schemes are $\cite[\mathrm{Section}\ 2]{Cor}$ and $\cite[\mathrm{Section}\ 2]{Ros13}$. Since the definition of jet schemes involves the Weil restriction, we will recall the definition and some properties of Weil restriction in Subsection 2.2 and then introduce the jet schemes in Subsection 2.3.

\subsection{The geometric version}
Before the discussion, we would like to mention that $f_{\mathcal{O}_{K}}$ in the statement of Theorem 1.0.3 is finite. Since $\mathrm{Pic}(\mathbb{P}_{\mathcal{O}_{K}}^{N})=\mathrm{Pic}(\mathcal{O}_{K})\times\mathbb{Z}=\mathbb{Z}$ which is generated by the twisting sheaf $\mathcal{O}(1)$ (see $\cite[(\mathrm{II},\ \mathrm{Ex.}\ 6.1,\ 6.16,\ 6.11.1\mathrm{A})]{Har77}$ and $\cite[\mathrm{Theorem}\ 40\ \mathrm{on}\ \mathrm{p.}\ 126]{Mat80}$), the $\mathcal{O}_{K}$-morphism $f_{\mathcal{O}_{K}}:\mathbb{P}_{\mathcal{O}_{K}}^{N}\rightarrow\mathbb{P}_{\mathcal{O}_{K}}^{N}$ is just determined by $N+1$ homogeneous polynomials in $\mathcal{O}_{K}[x_{0},\cdots,x_{N}]$. We know $f_{\mathcal{O}_{K}}$ will be finite if those polynomials are homogeneous of positive degree because it is affine under such situation (and it is always projective). But if those polynomials are constants, $f_{0}$ will be a constant map which contradicts to the second condition. So $f_{\mathcal{O}_{K}}$ is finite.

To prove Theorem 1.0.3, we need the Proposition below which gives a criterion for the DML property. It may be well-known to the experts, and some of its differnent forms have appeared in $\cite[\mathrm{Subsection}\ 3.1.3]{BGT16}$ and $\cite[\mathrm{Proposition}\ 4.2]{Xie14}$. But we would like to include a proof for completeness.

\begin{proposition}
Let $X$ be a quasi-projective variety over a field $K$ and let $f$ be a $K$-endomorphism of $X$. If for every integral closed subvariety $Y$ of $X$ of positive dimension and every point $x\in X(K)$, $\mathcal{O}_{f}(x)\cap Y$ is dense in $Y$ implies $Y$ is $f$-periodic, then $f$ satisfies the DML property.
\end{proposition}

\begin{prf}
Let $V$ be the closed subvariety of $X$ in the definition of the DML property. Denote $V_{0}=\overline{\mathcal{O}_{f}(x)\cap V}$ and let $Y_{1},\cdots,Y_{m}$ be the irreducible components of $V_{0}$. Then $Y_{1},\cdots,Y_{m}$ can be viewed as integral closed subvarieties of $X$ and $\mathcal{O}_{f}(x)\cap Y_{i}$ is dense in $Y_{i}$ for each $i=1,2,\cdots,m$. So by the assumption, $Y_{i}$ is $f$-periodic if it is not a single point.

We know $\{n\in\mathbb{N}|\ f^{n}(x)\in V(K)\}=\bigcup\limits_{i=1}^{m} \{n\in\mathbb{N}|\ f^{n}(x)\in Y_{i}(K)\}$. So it suffices to show that each $\{n\in\mathbb{N}|\ f^{n}(x)\in Y_{i}(K)\}$ is a finite union of arithmetic progressions. If $Y_{i}$ is not a single point, the result follows from it is $f$-periodic. Moreover, the result is trivial when $Y_{i}$ is a single point. So we are done.
\end{prf}

Now we will use Proposition 2.1.1 to show that it suffices to prove Theorem 2.1.2 below in order to prove Theorem 1.0.3.

\begin{theorem}
Let $K$ be a local field of positive characteristic. Let $k$ be the residue field of $\mathcal{O}_{K}$ which is a finite field. Let $f_{\mathcal{O}_{K}}:\mathbb{P}_{\mathcal{O}_{K}}^{N}\rightarrow\mathbb{P}_{\mathcal{O}_{K}}^{N}$ be a finite $\mathcal{O}_{K}$-morphism which satisfies:\\
1. $f_{\mathcal{O}_{K}}^{*}(\Omega_{\mathbb{P}_{\mathcal{O}_{K}}^{N}/\mathcal{O}_{K}})\rightarrow\Omega_{\mathbb{P}_{\mathcal{O}_{K}}^{N}/\mathcal{O}_{K}}$ is the zero map.\\
2. $f_{0}(=f_{\mathcal{O}_{K}}\times_{\mathcal{O}_{K}}k):\mathbb{P}_{k}^{N}\rightarrow\mathbb{P}_{k}^{N}$ is some $\mathrm{Frob}_{q}$ ($q$ is a power of the prime $\mathrm{char}k$).\\
Then the $K$-morphism $f:\mathbb{P}_{K}^{N}\rightarrow\mathbb{P}_{K}^{N}$ induced by $f_{\mathcal{O}_{K}}$ satisfies the property that for every geometrically integral closed subvariety $X\subseteq\mathbb{P}_{K}^{N}$ of positive dimension and every point $x\in\mathbb{P}^{N}(K)$, $X$ is $f$-periodic if $\mathcal{O}_{f}(x)\cap X$ is dense in $X$.
\end{theorem}

\begin{lemma}
Theorem 2.1.2 implies Theorem 1.0.3.
\end{lemma}

\begin{prf}
By Proposition 2.1.1, it suffices to verify that for every integral closed subvariety $X_{\overline{K}}\subseteq\mathbb{P}_{\overline{K}}^{N}$ of positive dimension and every point $x'\in\mathbb{P}^{N}(\overline{K})$, $X_{\overline{K}}$ is $f_{\overline{K}}$-periodic if $\mathcal{O}_{f_{\overline{K}}}(x')\cap X_{\overline{K}}$ is dense in $X_{\overline{K}}$ in the situation of Theorem 1.0.3.

Firstly, we know that both $X_{\overline{K}}$ and $x'$ come from corresponding objects defined over a finite extension of $K$. To be precise, there exists $K'$ which is a finite extension of $K$ such that we can find a closed subvariety $X_{K'}\subseteq\mathbb{P}_{K'}^{N}$ satisfies $X_{\overline{K}}=X_{K'}\times_{K'}\overline{K}$ and $x'\in\mathbb{P}^{N}(K')$. So $X_{K'}$ is a geometrically integral closed subvariety of $\mathbb{P}_{K'}^{N}$ of positive dimension. We will use Theorem 2.1.2 for the local field $K'$, in which $f_{\mathcal{O}_{K'}}$ is obtained from $f_{\mathcal{O}_{K}}$ by base-change. One can see that $f_{\mathcal{O}_{K'}}$ does satisfy the two conditions.

Since $\mathcal{O}_{f_{\overline{K}}}(x')\cap X_{\overline{K}}$ is dense in $X_{\overline{K}}$, we have $\mathcal{O}_{f_{K'}}(x')\cap X_{K'}$ is dense in $X_{K'}$ because the projection map $X_{\overline{K}}\rightarrow X_{K'}$ is surjective. Hence $X_{K'}$ if $f$-periodic in which $f:\mathbb{P}_{K'}^{N}\rightarrow\mathbb{P}_{K'}^{N}$ comes from $f_{\mathcal{O}_{K'}}$. As a result, we have $X_{\overline{K}}$ is $f_{\overline{K}}$-periodic by taking inverse image under the projection $\mathbb{P}_{\overline{K}}^{N}\rightarrow\mathbb{P}_{K'}^{N}$.
\end{prf}

So our task has turned into proving Theorem 2.1.2 from now on.

\begin{remark}
Notice that $K$ is isomorphic to $k((t))$ as a valuation field in the statement of Theorem 2.1.2. As a result, we can and do assume that $K=k((t))$ because it suffices to prove this case in order to prove the general one.
\end{remark}

\begin{remark}
Notice that $f^{m}$ satisfies the requirement in Theorem 2.1.2 for some positive integer $m$ implies that $f$ also satisfies since $X$ in there is an integral closed subvariety. Moreover, $f_{\mathcal{O}_{K}}^{m}$ is also finite and satisfies the first condition. Additionally, $f_{0}^{m}=\mathrm{Frob}_{q^m}$ in the second condition. Thus we may assume $q$ is a power of $\mathrm{|}k\mathrm{|}$ in Theorem 2.1.2 without loss of generality.
\end{remark}

\subsection{Weil restriction}
\begin{definition}
Let $S'\rightarrow S$ be a morphism of schemes. Let $X'$ be an $S'$-scheme. If the functor $\mathfrak{R}_{S'/S}(X'):(Sch/S)^{op}\rightarrow(Sets),T\mapsto\mathrm{Hom}_{S'}(T\times_{S}S',X')$ is representable by an $S$-scheme $X$, then we say that $X$ is the Weil restriction of $X'$ and also denote $X=\mathfrak{R}_{S'/S}(X')$.
\end{definition}

\begin{remark}
Keep the terminology in the definition, we can see that there is a functorial isomorphism $\mathrm{Hom}_{S}(T,\mathfrak{R}_{S'/S}(X'))\tilde{\rightarrow}\mathrm{Hom}_{S'}(T\times_{S}S',X')$.
\end{remark}

There is a criterion of the existence of Weil restriction scheme. See $\cite[7.6,\mathrm{Theorem}\ 4]{BLR90}$.

\begin{theorem}
Let $S'\rightarrow S$ be a morphism of schemes which is finite, flat and of finite presentation. Let $X'$ be an $S'$-scheme satisfies that for each $s\in S$ and each finite set of points $P\subseteq X'\times_{S}k(s)$, there is an open affine subscheme of $X'$ containing $P$. Then the Weil restriction of $X'$ exists.
\end{theorem}

\begin{remark}
One can check that the condition of Theorem 2.2.3 holds if $X'$ is a quasi-projective $S'$-scheme. Therefore, under the assumption that all schemes are noetherian, the Theorem says that the Weil restriction scheme exists if $S'\rightarrow S$ is finite flat and $X'$ is quasi-projective over $S'$.
\end{remark}

We know that the Weil restriction is naturally functorial in $X'$, namely, if $X_{1}'\rightarrow X_{2}'$ is an $S'$-morphism (and both Weil restriction schemes exist), then there is a natural $S$-morphism $\mathfrak{R}_{S'/S}(X_{1}')\rightarrow\mathfrak{R}_{S'/S}(X_{2}')$. We want to emphasize the fact that for an element which lies in $\mathrm{Hom}_{S}(T,\mathfrak{R}_{S'/S}(X_{1}'))=\mathrm{Hom}_{S'}(T\times_{S}S',X_{1}')$, it maps to the same element in $\mathrm{Hom}_{S}(T,\mathfrak{R}_{S'/S}(X_{2}'))=\mathrm{Hom}_{S'}(T\times_{S}S',X_{2}')$ by composing the above two morphisms, in which $T$ is an arbitrary $S$-scheme.

Moreover, this operation sends closed immersion to closed immersion. This is a special case (the case when Weil restriction scheme exists) of $\cite[7.6,\mathrm{Proposition}\ 2(\mathrm{ii})]{BLR90}$. We would like to state this result as a Proposition for completeness.

\begin{proposition}
Let $S'\rightarrow S$ be a morphism of schemes which is finite, flat and of finite presentation. Let $X_{1}'\rightarrow X_{2}'$ be a closed immersion of quasi-projective $S'$-schemes. Then the $S$-morphism $\mathfrak{R}_{S'/S}(X_{1}')\rightarrow\mathfrak{R}_{S'/S}(X_{2}')$ is also a closed immersion.
\end{proposition}

\subsection{Jet schemes}
For our purpose, we need to study the jet schemes over complete rings rather than jet schemes over finite type $k$-schemes. As a result, we will mainly follow the setting in $\cite[\mathrm{Section}\ 2]{Cor}$ in this subsection. But we should point out that in fact, many statements and proofs are just the same as their analogues in $\cite[\mathrm{Section}\ 2]{Ros13}$.

Denote $R=\mathcal{O}_{K}$ in which $K=k((t))$ as in Remark 2.1.4. Hence $\widehat{R^{\mathrm{sh}}}$ can be viewed as $\bar{k}[[t]]$. Denote $S=\mathrm{Spec}(\bar{k}[[t]])$ and let $S\widehat{\times}S=\mathrm{Spec}(\bar{k}[[x,y]])$ so that there exists a natural closed diagonal immersion $S\hookrightarrow S\widehat{\times}S$. Denote $S_{(n)}$ as the $n$th infinitesimal neighborhood of $S$ in $S\widehat{\times}S$, that is, $S_{(n)}=\mathrm{Spec}(\bar{k}[[x,y]]/((x-y)^{n+1}))$. Write $\pi_1,\pi_2:S\widehat{\times}S\rightarrow S$ for the first and second projection morphism and write $\pi_{1}^{S_{(n)}},\pi_{2}^{S_{(n)}}:S_{(n)}\rightarrow S$ for the induced morphisms. We view $S_{(n)}$ as a $S$-scheme via the \emph{first} projection. One can see that $S_{(n)}$ is finite flat over $S$ because $\bar{k}[[x,y]]/((x-y)^{n})$ is a free-of-rank-$n$ $\bar{k}[[t]]$-algebra. Moreover, all schemes (or rings) involved here are noetherian.

\begin{definition}
Let $W$ be a quasi-projective $S$-scheme. We define $\mathfrak{R}_{S_{(n)}/S}(\pi_{2}^{S_{(n)},*}W)$ as the $n$th jet scheme of $W$ over $S$ and denote by $J^{n}(W/S)$. Notice that the existence is guaranteed by Remark 2.2.4.
\end{definition}

We can see that the jet scheme is naturally functorial in $W$ and sends closed immersion to closed immersion as well as Weil restriction (Proposition 2.2.5). More precisely, an $S$-morphism $f:W\rightarrow W_{1}$ can induce an $S$-morphism $J^{n}(f):J^{n}(W/S)\rightarrow J^{n}(W_{1}/S)$ and $J^{n}(f)$ is a closed immersion if $f$ is. Just the same as we have mentioned in last subsection, composing $J^{n}(f)$ or $f_{n}:\pi_{2}^{S_{(n)},*}W\rightarrow\pi_{2}^{S_{(n)},*}W_{1}$ (the lifting of $f$) sends an element of $\mathrm{Hom}_{S}(T,J^{n}(W/S))=\mathrm{Hom}_{S_{(n)}}(T\times_{S}S_{(n)},\pi_{2}^{S_{(n)},*}W)$ to the same element in $\mathrm{Hom}_{S}(T,J^{n}(W_{1}/S))=\mathrm{Hom}_{S_{(n)}}(T\times_{S}S_{(n)},\pi_{2}^{S_{(n)},*}W_{1})$, in which $T$ is an arbitrary $S$-scheme.

For an $S$-scheme $T$, we may denote $T_{0}$ as the special fiber of $T$. We denote $S_{n}$ as the $n$th infinitesimal neighborhood of the closed point in $S$, that is, $S_{n}=\mathrm{Spec}(\bar{k}[[t]]/(t^{n+1}))$.

\begin{proposition}
Keep the terminology as above. There are canonical bijections $J^{n}(W/S)_{0}(\bar{k})=J^{n}(W/S)(S_{0})=W(S_{n})$.
\end{proposition}

\begin{prf}
By definition, $J^{n}(W/S)(S_{0})=\mathrm{Hom}_{S_{(n)}}(S_{0}\times_{S}S_{(n)},\pi_{2}^{S_{(n)},*}W)$. One can see that $S_{0}\times_{S}S_{(n)}=S_{n}$ as schemes. Moreover, compositing the structure morphism (projection) $S_{n}\rightarrow S_{(n)}$ with $\pi_{2}$, we just get the natural closed immersion $S_{n}\hookrightarrow S$. Hence $\mathrm{Hom}_{S_{(n)}}(S_{0}\times_{S}S_{(n)},\pi_{2}^{S_{(n)},*}W)=W(S_{n})$ and the result follows.
\end{prf}

We notice that there are canonical $S$-morphisms $\Lambda_{n,m}^{W}:J^{n}(W/S)\rightarrow J^{m}(W/S)$ induced by closed immersion $S_{(m)}\hookrightarrow S_{(n)}$ for any $n\geq m\geq0$. The reason is that we have $\mathrm{Hom}_{S_{(n)}}(T\times_{S}S_{(n)},\pi_{2}^{S_{(n)},*}W)=\mathfrak{R}_{S_{(n)}/S}(\pi_{2}^{S_{(n)},*}W)(T)$ and $\mathrm{Hom}_{S_{(m)}}(T\times_{S}S_{(m)},\pi_{2}^{S_{(m)},*}W)=\mathfrak{R}_{S_{(m)}/S}(\pi_{2}^{S_{(m)},*}W)(T)$ for every $S$-scheme $T$, and there is a natural restiction map $\mathrm{Hom}_{S_{(n)}}(T\times_{S}S_{(n)},\pi_{2}^{S_{(n)},*}W)\rightarrow\mathrm{Hom}_{S_{(m)}}(T\times_{S}S_{(m)},\pi_{2}^{S_{(m)},*}W)$. Again, as above, one can see that the map $\mathrm{Hom}_{S}(T,J^{n}(W/S))\rightarrow\mathrm{Hom}_{S}(T,J^{m}(W/S))$ given by composing $\Lambda_{n,m}^{W}$ is just the restriction map for arbitrary $S$-scheme $T$.

We will see that these morphisms have good properties when $W$ is smooth over $S$, or at least when $W=\mathbb{P}_{S}^{N}$.

\begin{proposition}
The morphism $\Lambda_{n,0}^{\mathbb{P}_{S}^{N}}$ is separated, of finite type and surjective. In particular, $J^{n}(\mathbb{P}_{S}^{N}/S)$ is a noetherian scheme for all $n\in\mathbb{N}$. 
\end{proposition}

\begin{prf}
To show $\Lambda_{n,0}^{\mathbb{P}_{S}^{N}}$ is separated and of finite type, it suffices to prove $J^{n}(\mathbb{P}_{S}^{N}/S)$ is a separated $S$-scheme of finite type. But this follows from $\cite[7.6,\mathrm{Proposition}\ 5(\mathrm{b})(\mathrm{e})]{BLR90}$ immediately.

To show $\Lambda_{n,0}^{\mathbb{P}_{S}^{N}}$ is surjective, it suffices to prove that it admits a section. Since $\mathbb{P}_{S}^{N}\times_{S}S_{(n)}=\mathbb{P}_{\mathbb{Z}}^{N}\times_{\mathbb{Z}}S_{(n)}=\pi_{2}^{S_{(n)},*}(\mathbb{P}_{S}^{N})$ as $S_{(n)}$-schemes, we have an $S_{(n)}$-isomorphism $i:\mathbb{P}_{S}^{N}\times_{S}S_{(n)}\tilde{\rightarrow}\pi_{2}^{S_{(n)},*}(\mathbb{P}_{S}^{N})$ such that $i\times_{S_{(n)}}S=\mathrm{id}_{\mathbb{P}_{S}^{N}}$. Then the element in $\mathrm{Hom}_{S}(\mathbb{P}_{S}^{N},J^{n}(\mathbb{P}_{S}^{N}/S))$ corresponding to $i$ is a section of $\Lambda_{n,0}^{\mathbb{P}_{S}^{N}}$.
\end{prf}

In fact, one can show that the morphism $\Lambda_{n,n-1}^{W}$ makes $J^{n}(W/S)$ into a $J^{n-1}(W/S)$-torsor under some vector bundle when $W$ is quasi-projective and smooth over $S$. See $\cite[\mathrm{Lemma}\ 2.2]{Cor}$.

At the end of this section, we mention that there is a canonical map $\lambda_{n}^{W}:W(S)\rightarrow J^{n}(W/S)(S)$ for any quasi-projective $S$-scheme $W$ given by the functoriality of jet schemes. More precisely, $\lambda_{n}^{W}(f):=J^{n}(f)$ (notice that $J^{n}(S/S)=S$). One can verify that this map satisfies $\Lambda_{n,0}^{W}\circ\lambda_{n}^{W}=\mathrm{id}$ by functorial properties. See $\cite[\mathrm{Lemma}\ 2.4]{Cor}$ or $\cite[\mathrm{Lemma}\ 2.5(\mathrm{a})]{Ros13}$.

\section{Critical schemes}
Recall that in order to prove Theorem 2.1.2, we have to show that for every geometrically integral closed subvariety $X\subseteq\mathbb{P}_{K}^{N}$ of positive dimension and every $x\in\mathbb{P}^{N}(K)$, $X$ is $f$-periodic if $\mathcal{O}_{f}(x)\cap X$ is dense in $X$ (in which $K=k((t))$ as we have assumed in Remark 2.1.4). We will lift the data from $R=\mathcal{O}_{K}$ and $K$ to $\widehat{R^{\mathrm{sh}}}$ and $L:=\widehat{K^{\mathrm{sh}}}=\bar{k}((t))$ at first in order to make use of the jet schemes introduced in last subsection. We will keep the terminology as above.

Denote $\mathcal{X}_{R}$ as the scheme-theoretic image of $X$ in $\mathbb{P}_{R}^{N}$. We can see that $\mathcal{X}_{R}$ is an integral closed subscheme in $\mathbb{P}_{R}^{N}$ and it is a model of $X$. Observe that $\mathcal{X}_{R}$ is flat over $\mathrm{Spec}(R)$ by $\cite[(\mathrm{III},(9.7))]{Har77}$. Denote $\mathcal{X}=\mathcal{X}_{R}\times_{R}S$ which is a closed subscheme in $\mathbb{P}_{S}^{N}$ and denote $X_{L}=X\times_{K}L$ as the generic fiber of $\mathcal{X}$ over $S$. Since $X$ is geometrically integral, $X_{L}$ is an integral closed subvariety of $\mathbb{P}_{L}^{N}$. Moreover, we have that $\mathcal{X}$ is integral since $X_{L},S$ are integral and $\mathcal{X}$ is flat over $S$ (see $\cite[4.3.1,\mathrm{Proposition}\ 3.8]{Liu02}$).

We lift the endomorphism $f_{R}:\mathbb{P}_{R}^{N}\rightarrow\mathbb{P}_{R}^{N}$ in the statement of Theorem 2.1.2 to $f_{S}:\mathbb{P}_{S}^{N}\rightarrow\mathbb{P}_{S}^{N}$ and $f_{L}:\mathbb{P}_{L}^{N}\rightarrow\mathbb{P}_{L}^{N}$. Notice that $f_{S}$ is also finite and satisfies the two conditions in Theorem 2.1.2, i.e. $f_{S}^{*}(\Omega_{\mathbb{P}_{S}^{N}/S})\rightarrow\Omega_{\mathbb{P}_{S}^{N}/S}$ is the zero map and $f_{S,0}:\mathbb{P}_{\bar{k}}^{N}\rightarrow\mathbb{P}_{\bar{k}}^{N}$ is $\mathrm{Frob}_{q}$.

Since both $\mathcal{X}$ and $\mathbb{P}_{S}^{N}$ are projective over $S$, we may consider the jet schemes $J^{n}(\mathcal{X}/S)$ and $J^{n}(\mathbb{P}_{S}^{N}/S)$ in the future. Moreover, $J^{n}(\mathcal{X}/S)$ is a closed subscheme of $J^{n}(\mathbb{P}_{S}^{N}/S)$. In addition, the jet schemes of $\mathbb{P}_{S}^{N}$ and the morphisms $\Lambda_{n,0}^{\mathbb{P}_{S}^{N}}$ have the good properties discussed in Proposition 2.3.3.

We will imitate some arguments in $\cite[\mathrm{Section}\ 3]{Ros13}$ in this section.

\subsection{Definition of Critical schemes}
Our key construction is the $S$-morphism $[f_{S}^{n}]^{\circ}:\mathbb{P}_{S}^{N}\rightarrow J^{n}(\mathbb{P}_{S}^{N}/S)$ which satisfies two conditions:
\begin{enumerate}
\item
$\Lambda_{n,0}^{\mathbb{P}_{S}^{N}}\circ[f_{S}^{n}]^{\circ}=f_{S}^{n}$
\item
$[f_{S}^{n}]^{\circ}\circ\Lambda_{n,0}^{\mathbb{P}_{S}^{N}}=J^{n}(f_{S})^{n}$
\end{enumerate}

Firstly, we have an $S_{(n)}$-isomorphism $i:\mathbb{P}_{S}^{N}\times_{S}S_{(n)}\tilde{\rightarrow}\pi_{2}^{S_{(n)},*}(\mathbb{P}_{S}^{N})$ such that $i\times_{S_{(n)}}S=\mathrm{id}_{\mathbb{P}_{S}^{N}}$ as in the proof of Proposition 2.3.3. Since $\mathrm{Hom}_{S_{(n)}}(\mathbb{P}_{S}^{N}\times_{S}S_{(n)},\pi_{2}^{S_{(n)},*}(\mathbb{P}_{S}^{N}))=\mathrm{Hom}_{S}(\mathbb{P}_{S}^{N},J^{n}(\mathbb{P}_{S}^{N}/S))$, we may define $[f_{S}^{n}]^{\circ}\in\mathrm{Hom}_{S}(\mathbb{P}_{S}^{N},J^{n}(\mathbb{P}_{S}^{N}/S))$ as the element which corresponding to $f_{S,n}^{n}\circ i\in\mathrm{Hom}_{S_{(n)}}(\mathbb{P}_{S}^{N}\times_{S}S_{(n)},\pi_{2}^{S_{(n)},*}(\mathbb{P}_{S}^{N}))$, in which $f_{S,n}:\pi_{2}^{S_{(n)},*}(\mathbb{P}_{S}^{N})\rightarrow\pi_{2}^{S_{(n)},*}(\mathbb{P}_{S}^{N})$ is the lifting of $f_{S}$.

\begin{lemma}
For each nonnegative integer $n$, the $S$-morphism $[f_{S}^{n}]^{\circ}$ satisfies the two conditions $\Lambda_{n,0}^{\mathbb{P}_{S}^{N}}\circ[f_{S}^{n}]^{\circ}=f_{S}^{n}$ and $[f_{S}^{n}]^{\circ}\circ\Lambda_{n,0}^{\mathbb{P}_{S}^{N}}=J^{n}(f_{S})^{n}$.
\end{lemma}

\begin{prf}
The reason of the first equality $\Lambda_{n,0}^{\mathbb{P}_{S}^{N}}\circ[f_{S}^{n}]^{\circ}=f_{S}^{n}$ is that $\Lambda_{n,0}^{\mathbb{P}_{S}^{N}}\circ[f_{S}^{n}]^{\circ}\in\mathrm{Hom}_{S}(\mathbb{P}_{S}^{N},\mathbb{P}_{S}^{N})$ is just the restriction of $f_{S,n}^{n}\circ i\in\mathrm{Hom}_{S_{(n)}}(\mathbb{P}_{S}^{N}\times_{S}S_{(n)},\pi_{2}^{S_{(n)},*}(\mathbb{P}_{S}^{N}))$, which is $f_{S}^{n}$. For the second equality, we need the hypothesis that $f_{S}^{*}(\Omega_{\mathbb{P}_{S}^{N}/S})\rightarrow\Omega_{\mathbb{P}_{S}^{N}/S}$ is the zero map.

Recall that $\Lambda_{n,0}^{\mathbb{P}_{S}^{N}}$ is the restriction of $j\in\mathrm{Hom}_{S_{(n)}}(J^{n}(\mathbb{P}_{S}^{N}/S)\times_{S}S_{(n)},\pi_{2}^{S_{(n)},*}(\mathbb{P}_{S}^{N}))$ to $S$ in which $j$ corresponds to the identity map in $\mathrm{Hom}_{S}(J^{n}(\mathbb{P}_{S}^{N}/S),J^{n}(\mathbb{P}_{S}^{N}/S))$. So in other words, $\Lambda_{n,0}^{\mathbb{P}_{S}^{N}}=j\times_{S_{(n)}}S$. By functorial properties, it suffices to prove that $f_{S,n}^{n}\circ i\circ(\Lambda_{n,0}^{\mathbb{P}_{S}^{N}}\times_{S}S_{(n)})=f_{S,n}^{n}\circ j\in\mathrm{Hom}_{S_{(n)}}(J^{n}(\mathbb{P}_{S}^{N}/S)\times_{S}S_{(n)},\pi_{2}^{S_{(n)},*}(\mathbb{P}_{S}^{N}))$. But since $f_{S}^{*}(\Omega_{\mathbb{P}_{S}^{N}/S})\rightarrow\Omega_{\mathbb{P}_{S}^{N}/S}$ is the zero map, we just have to prove $(i\circ(\Lambda_{n,0}^{\mathbb{P}_{S}^{N}}\times_{S}S_{(n)}))\times_{S_{(n)}}S=j\times_{S_{(n)}}S$ by $\cite[(\mathrm{III},5.1)]{SGA1}$. The last equality is true since $i\times_{S_{(n)}}S=\mathrm{id}_{\mathbb{P}_{S}^{N}}$ and $j\times_{S_{(n)}}S=\Lambda_{n,0}^{\mathbb{P}_{S}^{N}}$. Thus we are done.
\end{prf}

\begin{remark}
The $S$-morphisms $[f_{S}^{n}]^{\circ}$ are finite because $f_{S}$ is finite and $J^{n}(\mathbb{P}_{S}^{N}/S)$ is separated over $\mathbb{P}_{S}^{N}$ for each $n\in\mathbb{N}$.
\end{remark}

Now we can define the critical schemes.

\begin{definition}
We define $\mathrm{Crit}^{n}(\mathcal{X},\mathbb{P}_{S}^{N}):=[J^{n}(f_{S})^{n}]_{*}(J^{n}(\mathbb{P}_{S}^{N}/S))\cap J^{n}(\mathcal{X}/S)$, in which the first term stands for the scheme-theoretic image of $J^{n}(f_{S})^{n}:J^{n}(\mathbb{P}_{S}^{N}/S)\rightarrow J^{n}(\mathbb{P}_{S}^{N}/S)$ and $\cap$ stands for the scheme-theoretic intersection of those two closed subschemes of $J^{n}(\mathbb{P}_{S}^{N}/S)$.
\end{definition}

\begin{proposition}
The critical schemes $\mathrm{Crit}^{n}(\mathcal{X},\mathbb{P}_{S}^{N})$ are finite over $\mathbb{P}_{S}^{N}$.
\end{proposition}

\begin{prf}
It suffices to prove that $[J^{n}(f_{S})^{n}]_{*}(J^{n}(\mathbb{P}_{S}^{N}/S))$ is finite over $\mathbb{P}_{S}^{N}$ for each positive integer $n$. Since $[f_{S}^{n}]^{\circ}\circ\Lambda_{n,0}^{\mathbb{P}_{S}^{N}}=J^{n}(f_{S})^{n}$, we have $[J^{n}(f_{S})^{n}]_{*}(J^{n}(\mathbb{P}_{S}^{N}/S))=[f_{S}^{n}]^{\circ}_{*}([\Lambda_{n,0}^{\mathbb{P}_{S}^{N}}]_{*}(J^{n}(\mathbb{P}_{S}^{N}/S)))$. But we know that all schemes involved here are noetherian, $\Lambda_{n,0}^{\mathbb{P}_{S}^{N}}$ is surjective and $\mathbb{P}_{S}^{N}$ is reduced, so $[\Lambda_{n,0}^{\mathbb{P}_{S}^{N}}]_{*}(J^{n}(\mathbb{P}_{S}^{N}/S))=\mathbb{P}_{S}^{N}$ and hence $[J^{n}(f_{S})^{n}]_{*}(J^{n}(\mathbb{P}_{S}^{N}/S))=[f_{S}^{n}]^{\circ}_{*}(\mathbb{P}_{S}^{N})$.

By Proposition 2.3.3 and $\cite[(\mathrm{II},\mathrm{Ex.}\ 4.4)]{Har77}$, we have $[f_{S}^{n}]^{\circ}_{*}(\mathbb{P}_{S}^{N})$ is proper over $\mathbb{P}_{S}^{N}$ since $\Lambda_{n,0}^{\mathbb{P}_{S}^{N}}\circ[f_{S}^{n}]^{\circ}=f_{S}^{n}$. But the properness of $[f_{S}^{n}]^{\circ}$ implies that $\mathbb{P}_{S}^{N}$ maps onto $[f_{S}^{n}]^{\circ}_{*}(\mathbb{P}_{S}^{N})$. Hence $[f_{S}^{n}]^{\circ}_{*}(\mathbb{P}_{S}^{N})$ is quasi-finite over $\mathbb{P}_{S}^{N}$ because $f_{S}$ is finite. To sum up, $[J^{n}(f_{S})^{n}]_{*}(J^{n}(\mathbb{P}_{S}^{N}/S))=[f_{S}^{n}]^{\circ}_{*}(\mathbb{P}_{S}^{N})$ is finite over $\mathbb{P}_{S}^{N}$ and we are done.
\end{prf}

\begin{remark}
We can see that the natural morphism $J^{n}(\mathbb{P}_{S}^{N}/S)\rightarrow[J^{n}(f_{S})^{n}]_{*}(J^{n}(\mathbb{P}_{S}^{N}/S))$ induced by $J^{n}(f_{S})^{n}$ is surjective by the proof above.
\end{remark}

Notice that the natural morphism $[J^{n}(f_{S})^{n}]_{*}(J^{n}(\mathbb{P}_{S}^{N}/S))\rightarrow[J^{n-1}(f_{S})^{n-1}]_{*}(J^{n-1}(\mathbb{P}_{S}^{N}/S))$ and $\Lambda_{n,n-1}^{\mathcal{X}},\Lambda_{n,n-1}^{\mathbb{P}_{S}^{N}}$ form commutative diagrams with the closed immersions. Hence for each positive integer $n$, we get an $S$-morphism $\mathrm{Crit}^{n}(\mathcal{X},\mathbb{P}_{S}^{N})\rightarrow\mathrm{Crit}^{n-1}(\mathcal{X},\mathbb{P}_{S}^{N})$ by them. This morphism is finite since each $\mathrm{Crit}^{n}(\mathcal{X},\mathbb{P}_{S}^{N})$ is finite over $\mathbb{P}_{S}^{N}$. Thus we get a sequence of finite $S$-morphisms:

$$\cdots\rightarrow\mathrm{Crit}^{2}(\mathcal{X},\mathbb{P}_{S}^{N})\rightarrow\mathrm{Crit}^{1}(\mathcal{X},\mathbb{P}_{S}^{N})\rightarrow\mathrm{Crit}^{0}(\mathcal{X},\mathbb{P}_{S}^{N})=\mathcal{X}$$

We denote $\mathrm{Exc}^{n}(\mathcal{X},\mathbb{P}_{S}^{N})$ as the scheme-theoretic image of the morphism $\mathrm{Crit}^{n}(\mathcal{X},\mathbb{P}_{S}^{N})\rightarrow\mathcal{X}$, which is a closed subscheme of $\mathcal{X}$.

\subsection{Application}
We recall our mission. We have a geometrically integral closed subvariety $X\subseteq\mathbb{P}_{K}^{N}$ of positive dimension and a point $x\in\mathbb{P}^{N}(K)$. Our goal is to prove that $\mathcal{O}_{f}(x)\cap X$ is dense in $X$ implies $X$ is $f$-periodic. At this point, we would like to mention that we can almost forget $R=\mathcal{O}_{K}$ and $K$, and just think the question at the level of $S$ and $L$. To be precise, if there exists a point $x\in\mathbb{P}^{N}(K)$ such that $\mathcal{O}_{f}(x)\cap X$ is dense in $X$, then $\mathcal{O}_{f_{L}}(x)\cap X_{L}=\pi^{-1}(\mathcal{O}_{f}(x)\cap X)$ is dense in $X_{L}$ (in which we regard $x$ as a point in $\mathbb{P}^{N}(L)$) since the projection map $\pi:X_{L}\rightarrow X$ is open (see $\cite[\mathrm{Lemma}\ 29.23.4]{Sta22}$). Moreover, if $X_{L}$ is $f_{L}$-periodic, then $X$ will be $f$-periodic immediately by projection. \emph{So our mission turns into proving $X_{L}$ is $f_{L}$-periodic under the assumption that there exists a point $x\in\mathbb{P}^{N}(L)$ such that $\mathcal{O}_{f_{L}}(x)\cap X_{L}$ is dense in $X_{L}$ from now on.} In fact, we will prove that $f_{L}(X_{L})=X_{L}$ under this assumption later.

We will prove that $\mathrm{Exc}^{n}(\mathcal{X},\mathbb{P}_{S}^{N})=\mathcal{X}$ for every $n\in\mathbb{N}$ and deduce a lifting proposition (Proposition 3.2.2) under the assumption that $\mathcal{O}_{f_{L}}(x)\cap X_{L}$ is dense in $X_{L}$ in this subsection.

We may identify $\mathbb{P}_{S}^{N}(S)$ with $\mathbb{P}_{L}^{N}(L)$, then the subset $\mathcal{X}(S)$ corresponds to the subset $X_{L}(L)$ under this identification.

\begin{proposition}
If there exists a point $x\in\mathbb{P}^{N}(L)$ such that $\mathcal{O}_{f_{L}}(x)\cap X_{L}$ is dense in $X_{L}$, then $\mathrm{Exc}^{n}(\mathcal{X},\mathbb{P}_{S}^{N})=\mathcal{X}$ for every $n\in\mathbb{N}$.
\end{proposition}

\begin{prf}
Firstly, we can see that $f_{L}^{n}(\mathcal{O}_{f_{L}}(x))\cap X_{L}$ is dense in $X_{L}$ for each $n\in\mathbb{N}$ since $X_{L}$ is integral and of positive dimension. We may identify $f_{L}^{n}(\mathcal{O}_{f_{L}}(x))\cap X_{L}=f_{L}^{n}(\mathcal{O}_{f_{L}}(x))\cap X_{L}(L)\subseteq\mathbb{P}_{L}^{N}(L)$ with $f_{S}^{n}(\mathcal{O}_{f_{S}}(x))\cap\mathcal{X}(S)\subseteq\mathbb{P}_{S}^{N}(S)$.

We have $f_{S}^{n}(\mathcal{O}_{f_{S}}(x))\cap\mathcal{X}(S)=\Lambda_{n,0}^{\mathbb{P}_{S}^{N}}(\lambda_{n}^{\mathbb{P}_{S}^{N}}(f_{S}^{n}(\mathcal{O}_{f_{S}}(x))\cap\mathcal{X}(S)))$ by the statement at the end of Section 2. One can check $\lambda_{n}^{\mathbb{P}_{S}^{N}}(f_{S}^{n}(\mathcal{O}_{f_{S}}(x))\cap\mathcal{X}(S))\subseteq\lambda_{n}^{\mathbb{P}_{S}^{N}}(f_{S}^{n}(\mathcal{O}_{f_{S}}(x)))\cap\lambda_{n}^{\mathcal{X}}(\mathcal{X}(S))(\subseteq J^{n}(\mathbb{P}_{S}^{N}/S)(S))$. Furthermore, $\lambda_{n}^{\mathbb{P}_{S}^{N}}(f_{S}^{n}(\mathcal{O}_{f_{S}}(x)))\cap\lambda_{n}^{\mathcal{X}}(\mathcal{X}(S))\subseteq\mathrm{Crit}^{n}(\mathcal{X},\mathbb{P}_{S}^{N})(S)$ since the maps in $\lambda_{n}^{\mathbb{P}_{S}^{N}}(f_{S}^{n}(\mathcal{O}_{f_{S}}(x)))$ factor through $[J^{n}(f_{S})^{n}]_{*}(J^{n}(\mathbb{P}_{S}^{N}/S))$ and the maps in $\lambda_{n}^{\mathcal{X}}(\mathcal{X}(S))$ factor through $J^{n}(\mathcal{X}/S)$. As a result, we deduce $f_{S}^{n}(\mathcal{O}_{f_{S}}(x))\cap\mathcal{X}(S)\subseteq\Lambda_{n,0}^{\mathbb{P}_{S}^{N}}(\mathrm{Crit}^{n}(\mathcal{X},\mathbb{P}_{S}^{N})(S))$.

Now, since $\Lambda_{n,0}^{\mathbb{P}_{S}^{N}}(\mathrm{Crit}^{n}(\mathcal{X},\mathbb{P}_{S}^{N})(S))=\Lambda_{n,0}^{\mathcal{X}}(\mathrm{Crit}^{n}(\mathcal{X},\mathbb{P}_{S}^{N})(S))$ (as a subset of $\mathbb{P}_{S}^{N}(S)$), we can see that $\Lambda_{n,0}^{\mathbb{P}_{S}^{N}}(\mathrm{Crit}^{n}(\mathcal{X},\mathbb{P}_{S}^{N})(S))\subseteq\mathrm{Exc}^{n}(\mathcal{X},\mathbb{P}_{S}^{N})(S)$ by the definition of $\mathrm{Exc}^{n}(\mathcal{X},\mathbb{P}_{S}^{N})$. So going back to subsets contained in $\mathbb{P}_{L}^{N}(L)$, we get $f_{L}^{n}(\mathcal{O}_{f_{L}}(x))\cap X_{L}\subseteq\mathrm{Exc}^{n}(\mathcal{X},\mathbb{P}_{S}^{N})_{\eta}(L)\subseteq\mathrm{Exc}^{n}(\mathcal{X},\mathbb{P}_{S}^{N})_{\eta}$ in which $\eta$ is the generic point of $S$ and $\mathrm{Exc}^{n}(\mathcal{X},\mathbb{P}_{S}^{N})_{\eta}$ is the generic fiber of $\mathrm{Exc}^{n}(\mathcal{X},\mathbb{P}_{S}^{N})$.

However, since $\mathrm{Exc}^{n}(\mathcal{X},\mathbb{P}_{S}^{N})_{\eta}$ is a closed subscheme of $X_{L}$ and $f_{L}^{n}(\mathcal{O}_{f_{L}}(x))\cap X_{L}$ is dense in $X_{L}$, $\mathrm{Exc}^{n}(\mathcal{X},\mathbb{P}_{S}^{N})_{\eta}$ must be $X_{L}$ itself. As a result, $\mathrm{Exc}^{n}(\mathcal{X},\mathbb{P}_{S}^{N})$ must be $\mathcal{X}$ itself since $\mathcal{X}$ is integral.
\end{prf}

Now we can prove the main result of this subsection.

\begin{proposition}
Under the same assumption as Proposition 3.2.1, for each point $P\in\mathcal{X}_{0}(\bar{k})$, there exists a compatible sequence $\{P_{n}|\ P_{n}\in\mathcal{X}(S_{n})\cap f_{S}^{n}(\mathbb{P}_{S}^{N}(S_{n})),n\in\mathbb{N}\}$ such that $P_0$ is the natural lifting of $P$. Recall that for an $S$-scheme $T$, we denote $T_0$ as the special fiber of $T$. 
\end{proposition}

\begin{prf}
Firstly, we recall that there is a sequence of finite $\bar{k}$-morphisms such that each morphism $\mathrm{Crit}^{n}(\mathcal{X},\mathbb{P}_{S}^{N})_{0}\rightarrow\mathcal{X}_{0}$ is surjective (Proposition 3.2.1):

$$\cdots\rightarrow\mathrm{Crit}^{2}(\mathcal{X},\mathbb{P}_{S}^{N})_{0}\rightarrow\mathrm{Crit}^{1}(\mathcal{X},\mathbb{P}_{S}^{N})_{0}\rightarrow\mathrm{Crit}^{0}(\mathcal{X},\mathbb{P}_{S}^{N})_{0}=\mathcal{X}_{0}$$

So the maps $\mathrm{Crit}^{n}(\mathcal{X},\mathbb{P}_{S}^{N})_{0}(\bar{k})\rightarrow\mathcal{X}_{0}(\bar{k})$ must be surjective as well. Now we will consider the set $\mathrm{Crit}^{n}(\mathcal{X},\mathbb{P}_{S}^{N})_{0}(\bar{k})$.

We may observe that $\mathrm{Crit}^{n}(\mathcal{X},\mathbb{P}_{S}^{N})_{0}(\bar{k})=([J^{n}(f_{S})^{n}]_{*}(J^{n}(\mathbb{P}_{S}^{N}/S)))_{0}(\bar{k})\cap J^{n}(\mathcal{X}/S)_{0}(\bar{k})=\{Q_{n}\in J^{n}(\mathcal{X}/S)_{0}(\bar{k})|\ \exists\widetilde{Q_{n}}\in J^{n}(\mathbb{P}_{S}^{N}/S)_{0}(\bar{k}),\mathrm{s.t.}\ Q_{n}=J^{n}(f_{S})_{0}^{n}\circ\widetilde{Q_{n}}\}$ using Remark 3.1.5, in which $J^{n}(f_{S})_{0}$ is the lifting of $J^{n}(f_{S})$. By Proposition 2.3.2, we can identify $\mathrm{Crit}^{n}(\mathcal{X},\mathbb{P}_{S}^{N})_{0}(\bar{k})$ with the set $\{P_{n}\in\mathcal{X}(S_{n})|\ \exists\widetilde{P_{n}}\in\mathbb{P}_{S}^{N}(S_{n}),\mathrm{s.t.}\ P_{n}=f_{S}^{n}\circ{\widetilde{P_{n}}}\}=\mathcal{X}(S_{n})\cap f_{S}^{n}(\mathbb{P}_{S}^{N}(S_{n}))$. Thus we have proved that each $P\in\mathcal{X}_{0}(\bar{k})$ can lift to a $P_{n}\in\mathcal{X}(S_{n})\cap f_{S}^{n}(\mathbb{P}_{S}^{N}(S_{n}))$ for arbitrary $n\in\mathbb{N}$.

Now we have to show that we can choose the $P_{n}$ carefully in order to let them form a compatible sequence. Equivalently, we have to show that we can choose an $n$th preimage of $P\in\mathcal{X}_{0}$ in $\mathrm{Crit}^{n}(\mathcal{X},\mathbb{P}_{S}^{N})_{0}$ for each $n\in\mathbb{N}$ such that they are all compatible. But this follows from the quasi-finiteness of each connecting morphism $\mathrm{Crit}^{n}(\mathcal{X},\mathbb{P}_{S}^{N})_{0}\rightarrow\mathrm{Crit}^{n-1}(\mathcal{X},\mathbb{P}_{S}^{N})_{0}$ and the surjectiveness of each $\mathrm{Crit}^{n}(\mathcal{X},\mathbb{P}_{S}^{N})_{0}\rightarrow\mathcal{X}_{0}$ immediately.
\end{prf}

\section{The proof of Theorem 1.0.3}
We will prove that $f_{L}(X_{L})=X_{L}$ in this Section and thus finish the proof. We must make use of the condition that $f_{S,0}:\mathbb{P}_{\bar{k}}^{N}\rightarrow\mathbb{P}_{\bar{k}}^{N}$ is $\mathrm{Frob}_{q}$. We absorb the ideas used in $\cite[\mathrm{Subsection}\ 4.5]{Xie18}$.

Notice we have assumed that $q$ is a power of $\mathrm{|}k\mathrm{|}$ in Remark 2.1.5, so we may denote $\sigma=\mathrm{Frob}_{q}\in\mathrm{Gal}(\bar{k}/k)$ (one should distinguish the element in Galois group and the endomorphism of $\mathbb{P}_{\bar{k}}^{N}$ although both of them can be written as $\mathrm{Frob}_{q}$). Since $\mathrm{Gal}(K^{\mathrm{sh}}/K)=\mathrm{Gal}(\bar{k}/k)$, it induces an isomorphism of $L=\bar{k}((t))$ which fixes elements in $K=k((t))$ (just the $\mathrm{Frob}_{q}$ acting on the coefficients). As a result, it induces a map $\mathbb{P}_{S}^{N}(S)\rightarrow\mathbb{P}_{S}^{N}(S)$ (we may identify $\mathbb{P}_{S}^{N}(S)$ with $\mathbb{P}_{L}^{N}(L)$) and hence a map $\mathbb{P}_{S}^{N}(S_{n})\rightarrow\mathbb{P}_{S}^{N}(S_{n})$ for each $n\in\mathbb{N}$ as well. Abusing notation, we may call all these maps by $\sigma$. 

Since $f_{S}$ comes from $f_{R}$, we can see that the operation $\sigma:\mathbb{P}_{S}^{N}(S)\rightarrow\mathbb{P}_{S}^{N}(S)$ (or $\mathbb{P}_{S}^{N}(S_{n})\rightarrow\mathbb{P}_{S}^{N}(S_{n})$) commutes with composing $f_{S}$. 

\begin{lemma}
Let $n\in\mathbb{N}$. For an element $P_{n}\in f_{S}^{n}(\mathbb{P}_{S}^{N}(S_{n}))$, we have $f_{S}\circ P_{n}=\sigma(P_{n})$.
\end{lemma}

\begin{prf}
Suppose that $P_{n}=f_{S}^{n}\circ Q_{n}$ for a $Q_{n}\in\mathbb{P}_{S}^{N}(S_{n})$. We only have to prove that $f_{S}^{n}\circ(f_{S}\circ Q_{n})=f_{S}^{n}\circ(\sigma(Q_{n}))$ because $f_{S}\circ\sigma=\sigma\circ f_{S}$. Since $f_{S}^{*}(\Omega_{\mathbb{P}_{S}^{N}/S})\rightarrow\Omega_{\mathbb{P}_{S}^{N}/S}$ is the zero map, we just have to show that $f_{S}\circ Q_{n}$ and $\sigma(Q_{n})$ reduce to the same element in $\mathbb{P}_{\bar{k}}^{N}(\bar{k})$ by $\cite[(\mathrm{III},5.1)]{SGA1}$. But this follows from $f_{S,0}=\mathrm{Frob}_{q}$.
\end{prf}

Now for each point $P\in\mathcal{X}_{0}(\bar{k})$, we can use the compatible sequence $\{P_{n}\}$ in Proposition 3.2.2 to construct an element $\widetilde{P}\in\mathcal{X}(S)(\subseteq\mathbb{P}_{S}^{N}(S))$ which satisfies $f_{S}\circ\widetilde{P}=\sigma(\widetilde{P})$ by Lemma 4.0.1. We may identify $\mathcal{X}(S)$ with $X_{L}(L)$ and define $\mathscr{P}=\{\widetilde{P}|\ P\in\mathcal{X}_{0}(\bar{k})\}\subseteq X_{L}(L)(\subseteq\mathbb{P}_{L}^{N}(L))$. Using $\cite[\mathrm{Lemma}\ 33.19.2]{Sta22}$, we can see that $\mathscr{P}$ is a dense subset in $X_{L}$ since the closed points are dense in $\mathcal{X}_{0}$. Regard $\mathscr{P}$ as a dense subset of prime ideals in $X_{L}\subseteq\mathbb{P}_{L}^{N}$. For each prime ideal $\mathfrak{p}\in\mathscr{P}$, we have $f_{L}(\mathfrak{p})=\sigma(\mathfrak{p})$ because the map in $\mathbb{P}_{L}^{N}(L)$ which corresponding to $\mathfrak{p}$ satisfies the same requirement.

Since $f_{L}$ is also a finite morphism, we deduce that $f_{L}(X_{L})=\overline{f_{L}(\mathscr{P})}=\overline{\sigma(\mathscr{P})}=\sigma(X_{L})$. Here $\sigma$ acts on closed subsets in $\mathbb{P}_{L}^{N}$ by acting on the coefficients of the defining equations. But $X_{L}$ comes from $X$ which is defined over $K$ and $\sigma$ fixes elements in $K$, so we have $f_{L}(X_{L})=\sigma(X_{L})=X_{L}$ and hence finish the proof of Theorem 2.1.2. As a result, we have proved Theorem 1.0.3 by taking Lemma 2.1.3 into account.

\section{The proof of Theorem 1.0.2}
We will prove Theorem 1.0.2 and a generalized version in this Section. But firstly, we would like to propose a corollary of Theorem 1.0.3 and provide some examples.

\subsection{A corollary and some examples}
\begin{corollary}
Let $K$ be the function field of a variety $V$ over $\overline{\mathbb{F}_{p}}$. Let $f_{V}:\mathbb{P}_{V}^{N}\rightarrow\mathbb{P}_{V}^{N}$ be a $V$-morphism which satisfies:\\
1. $f_{V}^{*}(\Omega_{\mathbb{P}_{V}^{N}/V})\rightarrow\Omega_{\mathbb{P}_{V}^{N}/V}$ is the zero map.\\
2. There exists a nonsingular closed point $u\in V$, such that $f_{u}(=f_{V}\times_{V}\{u\}):\mathbb{P}_{\overline{\mathbb{F}_{p}}}^{N}\rightarrow\mathbb{P}_{\overline{\mathbb{F}_{p}}}^{N}$ is some $\mathrm{Frob}_{q}$ ($q$ is a power of the prime $p$).\\
Then the $\overline{K}$-morphism $f_{\overline{K}}:\mathbb{P}_{\overline{K}}^{N}\rightarrow\mathbb{P}_{\overline{K}}^{N}$ induced by $f_{V}$ satisfies the DML property.
\end{corollary}

\begin{prf}
Firstly, we may substitute $V$ by a standard smooth open affine neighborhood of $u$ in $V$. Then we may choose an appropriate positive integer $m$ such that $V,f_{V}$ and $u$ come from $\overline{V},f_{\overline{V}}$ and $\overline{u}$ in which $\overline{V}$ is a smooth variety over $\mathbb{F}_{p^m},\overline{u}\in\overline{V}(\mathbb{F}_{p^m})$ and $f_{\overline{V}}$ is a $\overline{V}$-endomorphism of $\mathbb{P}_{\overline{V}}^{N}$ that satisfies the analogues of the two conditions, i.e. $f_{\overline{V}}^{*}(\Omega_{\mathbb{P}_{\overline{V}}^{N}/\overline{V}})\rightarrow\Omega_{\mathbb{P}_{\overline{V}}^{N}/\overline{V}}$ is the zero map and $f_{\overline{u}}:\mathbb{P}_{\mathbb{F}_{p^m}}^{N}\rightarrow\mathbb{P}_{\mathbb{F}_{p^m}}^{N}$ is $\mathrm{Frob}_{q}$. We denote $K'$ as the function field of $\overline{V}$ which is a finitely generated field over $\mathbb{F}_{p^m}$, and denote $f_{K'}$ as the endomorphism of $\mathbb{P}_{K'}^{N}$ induced by $f_{\overline{V}}$ on the generic fiber.

In order to make use of Theorem 1.0.3, we need to construct a morphism $\mathrm{Spec}(\mathcal{O}_{L})\rightarrow\overline{V}$ in which $L$ is a local field of characteristic $p$. This morphism should send the generic point in $\mathrm{Spec}(\mathcal{O}_{L})$ to the generic point in $\overline{V}$ and send the special point in $\mathrm{Spec}(\mathcal{O}_{L})$ to $\overline{u}$. By $\cite[(\mathrm{II},4.4)]{Har77}$, our mission turns into constructing a local field $L$ of characteristic $p$ that containing $K'$ such that $\mathcal{O}_{L}$ dominates $\mathcal{O}_{\overline{V},\overline{u}}$.

Denote $d=\mathrm{dim}\overline{V}=\mathrm{dim}\mathcal{O}_{\overline{V},\overline{u}}=\mathrm{tr.deg.}$ $K'/\mathbb{F}_{p^m}$. Since $\mathcal{O}_{\overline{V},\overline{u}}$ is a regular local ring, we may let $x_{1},\cdots,x_{d}$ be the local parameters in its maximal ideal. One can see that $x_{1},\cdots,x_{d}$ are algebraically independent over $\mathbb{F}_{p^m}$ which is the residue field of $\mathcal{O}_{\overline{V},\overline{u}}$, so $K'$ is a finite extension of $\mathbb{F}_{p^m}(x_{1},\cdots,x_{d})$. But tr.deg. $\mathbb{F}_{p^m}((t))/\mathbb{F}_{p^m}=\infty$ since the field of Laurant series is uncountable, so we may choose $g_{1},\cdots,g_{d}\in t+t^{2}\cdot\mathbb{F}_{p^m}[[t]]$ such that they are algebraically independent over $\mathbb{F}_{p^m}$. As a result, we obtain a homomorphism $\mathcal{O}_{\overline{V},\overline{u}}\rightarrow\mathbb{F}_{p^m}[[t]]$ by composing $\mathcal{O}_{\overline{V},\overline{u}}\hookrightarrow\widehat{\mathcal{O}_{\overline{V},\overline{u}}}=\mathbb{F}_{p^m}[[x_{1},\cdots,x_{d}]]\rightarrow\mathbb{F}_{p^m}[[t]]$ in which the latter map sends $x_{i}$ to $g_{i}$ for each $i=1,\cdots,d$. Since each element of $\mathcal{O}_{\overline{V},\overline{u}}$ is algebraic over $\mathbb{F}_{p^m}(x_{1},\cdots,x_{d})$ and $g_{1},\cdots,g_{d}$ are algebraically independent over $\mathbb{F}_{p^m}$, we can see that the homomorphism $\mathcal{O}_{\overline{V},\overline{u}}\rightarrow\mathbb{F}_{p^m}[[t]]$ is injective. Moreover, it sends the elements in the maximal ideal of $\mathcal{O}_{\overline{V},\overline{u}}$ to elements in $t\cdot\mathbb{F}_{p^m}[[t]]$ by definition. Thus we can just take $L=\mathbb{F}_{p^m}((t))$ which contains $\mathrm{Frac}(\mathcal{O}_{\overline{V},\overline{u}})=K'$ and get the morphism $\mathrm{Spec}(\mathcal{O}_{L})\rightarrow\overline{V}$.

Once we obtain the desired morphism $\mathrm{Spec}(\mathcal{O}_{L})\rightarrow\overline{V}$, we can use Theorem 1.0.3 to conclude that the endomorphism $f_{\overline{L}}$ of $\mathbb{P}_{\overline{L}}^{N}$ induced by $f_{\mathcal{O}_{L}}$ (or by $f_{\overline{V}}$ more essensially) satisfies the DML property. As a result, $f_{\overline{K}}$ also satisfies the DML property since $\overline{K}=\overline{K'}\subseteq\overline{L}$ and $f_{\overline{L}}$ is also induced by $f_{\overline{K'}}$. Thus we have finished the proof.
\end{prf}

\begin{remark}
If $V$ is a nonsingular variety over $\overline{\mathbb{F}_{p}}$ and $u\in V(\overline{\mathbb{F}_{p}})$, then we may construct a morphism $\mathrm{Spec}(\mathcal{O}_{L})\rightarrow V$ sending the generic point in $\mathrm{Spec}(\mathcal{O}_{L})$ to the generic point in $V$ and sending the special point in $\mathrm{Spec}(\mathcal{O}_{L})$ to $u$ in which $L=\overline{\mathbb{F}_{p}}((t))$. The proof is just the same as above.
\end{remark}

Now we provide two direct applications of our results.

\begin{example}
Let $K$ be a local field of characteristic $p>0$. Let $k=\overline{\mathbb{F}_{p}}\cap K$ be the coefficient field of $K$ (which is isomorphic to the finite residue field of $\mathcal{O}_K$). Then the endomorphism $f:\mathbb{P}_{\overline{K}}^{N}\rightarrow\mathbb{P}_{\overline{K}}^{N}$,
$$\begin{bmatrix} x_{0} \\ x_{1} \\ \vdots \\ x_{N} \end{bmatrix} \mapsto \begin{bmatrix} \sum\limits_{i=0}^{N}a_{0i}x_{i}^{q}+g_{0}(x_{0}^{p},\cdots,x_{N}^{p}) \\ \sum\limits_{i=0}^{N}a_{1i}x_{i}^{q}+g_{1}(x_{0}^{p},\cdots,x_{N}^{p}) \\ \vdots \\ \sum\limits_{i=0}^{N}a_{Ni}x_{i}^{q}+g_{N}(x_{0}^{p},\cdots,x_{N}^{p}) \end{bmatrix}$$
satisfies the DML property, in which $A=(a_{ij})_{(N+1)\times(N+1)}\in GL_{N+1}(k)$, $q$ is a power of $p$ and $g_{0},\cdots,g_{N}\in\mathfrak{m}_{K}[x_{0},\cdots,x_{N}]$ are homogeneous polynomials of degree $=\frac{q}{p}$.
\end{example}

\begin{prf}
Taking Remark 1.0.4 into account, this is just a special case of Theorem 1.0.3.
\end{prf}

\begin{example}
Let $K=\overline{\mathbb{F}_{p}}(t_{1},\cdots,t_{d})$. Then the endomorphism $f:\mathbb{P}_{\overline{K}}^{N}\rightarrow\mathbb{P}_{\overline{K}}^{N}$,
$$\begin{bmatrix} x_{0} \\ x_{1} \\ \vdots \\ x_{N} \end{bmatrix} \mapsto \begin{bmatrix} \sum\limits_{i=0}^{N}a_{0i}x_{i}^{q}+g_{0}(x_{0}^{p},\cdots,x_{N}^{p}) \\ \sum\limits_{i=0}^{N}a_{1i}x_{i}^{q}+g_{1}(x_{0}^{p},\cdots,x_{N}^{p}) \\ \vdots \\ \sum\limits_{i=0}^{N}a_{Ni}x_{i}^{q}+g_{N}(x_{0}^{p},\cdots,x_{N}^{p}) \end{bmatrix}$$
satisfies the DML property, in which $A=(a_{ij})_{(N+1)\times(N+1)}\in GL_{N+1}(\overline{\mathbb{F}_{p}})$, $q$ is a power of $p$ and $g_{0},\cdots,g_{N}\in\overline{\mathbb{F}_{p}}[t_{1},\cdots,t_{d}][x_{0},\cdots,x_{N}]$ are homogeneous polynomials of degree $=\frac{q}{p}$ such that every coefficient in $\overline{\mathbb{F}_{p}}[t_{1},\cdots,t_{d}]$ has zero constant term.
\end{example}

\begin{prf}
This is a consequence of Example 5.1.3 above.
\end{prf}

\subsection{The proof of Theorem 1.0.2 and a generalization}
Now we are going to prove Theorem 1.0.2. Firstly, we will prove a technical lemma in order to show that the matrix $A$ in the statement of Theorem 1.0.2 can be assumed to be the identity matrix without loss of generality. Then we will use an embedding argument to finish the proof.

\begin{lemma}
Let $K$ be an algebraically closed field of characteristic $p>0$. Let $q$ be a power of $p$. Then for every matrix $A\in GL_{N}(K)$, there exists a matrix $B\in GL_{N}(K)$ such that $B^{(q)}=AB$ in which $B^{(q)}$ is an abbreviation of $\mathrm{Frob}_{q}(B)$. In other words, we require each element of $AB$ to be the $q$th power of the corresponding element of $B$.
\end{lemma}

\begin{prf}
Equivalently, we have to prove that there are linearly independent vectors $\{\beta_{1},\cdots,\beta_{N}\}\subseteq K^{N}$ such that $\beta_{i}^{(q)}=A\beta_{i}$ for each $i$, in which $\beta^{(q)}=\mathrm{Frob}_{q}(\beta)$ has the same meaning as in the statement of this lemma. Since $A\in GL_{N}(K)$ and $K$ is an algebraically closed field of characteristic $p$, one can see that the equation $\beta^{(q)}=A\beta$ has exactly $q^{N}$ different solutions in $K^{N}$ by intersection theory (notice that it has only finitely many solutions and each solution has multiplicity 1).

Suppose that $\{\beta_{1},\cdots,\beta_{m}\}$ is a maximal linear indepent subset among those $q^{N}$ solutions. Then each solution $\beta$ has the form $\sum\limits_{i=1}^{m} c_{i}\beta_{i}$ in which each $c_{i}\in K$. Now $\beta^{(q)}=A\beta$ can be read as $\sum\limits_{i=1}^{m} c_{i}^{q}\beta_{i}^{(q)}=A\sum\limits_{i=1}^{m} c_{i}\beta_{i}$. So $A\sum\limits_{i=1}^{m} c_{i}^{q}\beta_{i}=A\sum\limits_{i=1}^{m} c_{i}\beta_{i}$ and hence $\sum\limits_{i=1}^{m} c_{i}^{q}\beta_{i}=\sum\limits_{i=1}^{m} c_{i}\beta_{i}$ since $A$ is invertible. Now the linear independence of $\{\beta_{1},\cdots,\beta_{m}\}$ implies that each $c_{i}\in\mathbb{F}_{q}$, so the equation has at most $q^{m}$ solutions. Thus $m\geq N$ and we are done.
\end{prf}

\begin{corollary}
In order to prove Theorem 1.0.2, we may assume that the matrix $A$ in there is the identity matrix without loss of generality.
\end{corollary}

\begin{prf}
By Lemma 5.2.1, we can find $B\in GL_{N+1}(K)$ such that $B^{(q)}=A^{-1}B$. Then $B$ lies in $GL_{N+1}(\mathcal{O}_{K})$ automatically since $A^{-1}\in GL_{N+1}(\mathcal{O}_{K})$. Denote $\sigma$ as the automorphism of $\mathbb{P}_{K}^{N}$ corresponding to $B$. Then one can see that $\sigma^{-1}f\sigma$ has the same form as $f$ and the matrix becomes the identity matrix. Thus we have finished the proof since $f$ will satisfy the DML property if $\sigma^{-1}f\sigma$ satisfies.
\end{prf}

Next, we shall use an embedding argument. We will denote $F$ as $\widehat{\overline{\mathbb{F}_{p}((t))}}$ and denote $L$ as $\overline{\mathbb{F}_{p}}((t))$ from now on. Notice that $L\subseteq F$.

\begin{lemma}
Let $a_{1},\cdots,a_{n}$ be $n$ elements in the maximal ideal $\mathfrak{m}_{F}\subseteq\mathcal{O}_{F}$. Then one can find a local field $K_{0}$ and an embedding $i:\mathbb{F}_{p}(a_{1},\cdots,a_{n})\hookrightarrow K_{0}$, such that $i(a_{1}),\cdots,i(a_{n})$ all lie in $\mathfrak{m}_{K_{0}}\subseteq\mathcal{O}_{K_{0}}$.
\end{lemma}

\begin{prf}
Firstly, using $\cite[(\mathrm{I},4.7\mathrm{A},4.8\mathrm{A})]{Har77}$, we may assume that $a_{1},\cdots,a_{m}$ is a set of separating transcendence base of $\mathbb{F}_{p}(a_{1},\cdots,a_{n})/\mathbb{F}_{p}$ without loss of generality. That is to say, $\{a_{1},\cdots,a_{m}\}$ are algebraically independent over $\mathbb{F}_{p}$ and $\mathbb{F}_{p}(a_{1},\cdots,a_{n})/\mathbb{F}_{p}(a_{1},\cdots,a_{m})$ is a finite separable extension. As a result, we can find an element $\theta\in\mathbb{F}_{p}(a_{1},\cdots,a_{n})$ such that $\mathbb{F}_{p}(a_{1},\cdots,a_{n})=\mathbb{F}_{p}(a_{1},\cdots,a_{m})(\theta)$. 

Let $\theta^{d}+f_{d-1}\theta^{d-1}+\cdots+f_{0}=0$ be the minimal polynomial of $\theta$ over $\mathbb{F}_{p}(a_{1},\cdots,a_{m})$. Write $a_{j}=\sum\limits_{k=0}^{d-1} g_{jk}\theta^{k}$ in which all $g_{jk}\in\mathbb{F}_{p}(a_{1},\cdots,a_{m})$ for each $j=1,2,\cdots,n$. Since $F_{0}=\overline{\mathbb{F}_{p}((t))}$ is dense in $F$ and local fields are uncountable, we can find $a_{1}',\cdots,a_{m}'\in F_{0}$ such that they are algebraically independent over $\mathbb{F}_{p}$ and each $a_{i}'$ is very very close to $a_{i}$. As a result, we can find a root $\theta'\in F_{0}$ of $x^{d}+f_{d-1}(a_{1}',\cdots,a_{m}')x^{d-1}+\cdots+f_{0}(a_{1}',\cdots,a_{m}')$ which is sufficiently close to $\theta$. Now let $K_{0}=\mathbb{F}_{p}((t))(a_{1}',\cdots,a_{m}')(\theta')$ which is a local field, then we get an embedding $i:\mathbb{F}_{p}(a_{1},\cdots,a_{n})\hookrightarrow K_{0}$.

Since $i(a_{j})=\sum\limits_{k=0}^{d-1} g_{jk}(a_{1}',\cdots,a_{m}')\theta'^{k}$ is very close to $a_{j}$ for each $j=1,2,\cdots,n$ and every $a_{j}$ has absolute value less than 1, we know that every $i(a_{j})$ also has absolute value less than 1. In other words, $i(a_{1}),\cdots,i(a_{n})$ all lie in $\mathfrak{m}_{K_{0}}$. Thus we are done.
\end{prf}

\begin{lemma}
Let $A=\overline{\mathbb{F}_{p}}[x_{1},\cdots,x_{n}]/\mathfrak{p}$ be an integral domain. If every polynomial in $\mathfrak{p}$ has zero constant term, then there exists an embedding $A\hookrightarrow\mathcal{O}_{L}$ sending each $x_{i}$ into $\mathfrak{m}_{L}$.
\end{lemma}

\begin{prf}
Denote $V=\mathrm{Spec}(A)$ which is a variety over $\overline{\mathbb{F}_{p}}$ and denote $u=(0,\cdots,0)\in V(\overline{\mathbb{F}_{p}})$. Using $\cite[\mathrm{Theorem}\ 3.1]{dJ96}$, we can find a nonsingular variety $V'$ over $\overline{\mathbb{F}_{p}}$ such that there exists a dominant proper $\overline{\mathbb{F}_{p}}$-morphism $\phi:V'\rightarrow V$. So $\phi$ sends the generic point of $V'$ to the generic point of $V$ and there exists a point $u'\in V'(\overline{\mathbb{F}_{p}})$ which maps to $u$. Combining with Remark 5.1.2, we get a morphism $\mathrm{Spec}(\mathcal{O}_{L})\rightarrow V$ sending the generic point in $\mathrm{Spec}(\mathcal{O}_{L})$ to the generic point in $V$ and sending the special point in $\mathrm{Spec}(\mathcal{O}_{L})$ to $u$. Thus we are done.
\end{prf}

Combining the two lemmas above, we deduce the Proposition below.

\begin{proposition}
Let $K$ be a complete algebraically closed non-archimedian valuation field of characteristic $p>0$. Let $a_{1},\cdots,a_{n}$ be $n$ elements in the maximal ideal $\mathfrak{m}_{K}\subseteq\mathcal{O}_{K}$. Then one can find a local field $K_{0}$ and an embedding $i:\mathbb{F}_{p}(a_{1},\cdots,a_{n})\hookrightarrow K_{0}$, such that $i(a_{1}),\cdots,i(a_{n})$ all lie in $\mathfrak{m}_{K_{0}}\subseteq\mathcal{O}_{K_{0}}$.
\end{proposition}

Now we can prove Theorem 1.0.2.

\proof[Proof of Theorem 1.0.2]
By Corollary 5.2.2, we may assume that the matrix $A$ is the identity matrix. Let $a_{1},\cdots,a_{n}\in\mathfrak{m}_{K}$ be all coefficients of $g_{0},\cdots,g_{N}$. Using Proposition 5.2.5, we may embed $\mathbb{F}_{p}(a_{1},\cdots,a_{n})$ into a local field $K_{0}$ such that each $a_{i}$ goes into $\mathfrak{m}_{K_{0}}$ under this embedding. We want to show that for any point $x\in\mathbb{P}^{N}(K)$ and any closed subvariety $V\subseteq\mathbb{P}_{K}^{N}$, the set $\{n\in\mathbb{N}|\ f^{n}(x)\in V\}$ is a finite union of arithmetic progressions.

Let $x_{0},\cdots,x_{N}$ be the coefficients of $x$ and let $y_{1},\cdots,y_{m}$ be the coefficients of the defining equations of $V$. Since $K_{0}$ has infinite transcendence degree over $\mathbb{F}_{p}$, we can extend the embedding $\mathbb{F}_{p}(a_{1},\cdots,a_{n})\hookrightarrow K_{0}$ to an embedding $\mathbb{F}_{p}(a_{1},\cdots,a_{n},x_{0},\cdots,x_{N},y_{1},\cdots,y_{m})\hookrightarrow K_{0}'$ in which $K_{0}'$ is a finite extension of $K_{0}$. So we have descended all the data to the local field $K_{0}'$ and thus the result follows from Example 5.1.3.
\endproof

We would like to mention that in fact we have shown a more general statement by the embedding argument above. But since it is not as succinct as Theorem 1.0.2, we will write it down as a Proposition below. 

\begin{proposition}
Let $K$ be an algebraically closed field of characteristic $p>0$. Let $f:\mathbb{P}_{K}^{N}\rightarrow\mathbb{P}_{K}^{N}$,
$$\begin{bmatrix} x_{0} \\ x_{1} \\ \vdots \\ x_{N} \end{bmatrix} \mapsto \begin{bmatrix} \sum\limits_{i=0}^{N}a_{0i}x_{i}^{q}+g_{0}(x_{0}^{p},\cdots,x_{N}^{p}) \\ \sum\limits_{i=0}^{N}a_{1i}x_{i}^{q}+g_{1}(x_{0}^{p},\cdots,x_{N}^{p}) \\ \vdots \\ \sum\limits_{i=0}^{N}a_{Ni}x_{i}^{q}+g_{N}(x_{0}^{p},\cdots,x_{N}^{p}) \end{bmatrix}$$
be an endomorphism in which $A=(a_{ij})_{(N+1)\times(N+1)}\in GL_{N+1}(\overline{\mathbb{F}_{p}})$, $q$ is a power of $p$ and $g_{0},\cdots,g_{N}\in K[x_{0},\cdots,x_{N}]$ are homogeneous polynomials of degree $=\frac{q}{p}$. Let $a_{1},\cdots,a_{n}$ be all of the coefficients of $g_{0},\cdots,g_{N}$. If $g(a_{1},\cdots,a_{n})\neq0$ for every polynomial $g\in\overline{\mathbb{F}_{p}}[y_{1},\cdots,y_{n}]$ with nonzero constant term, then $f$ satisfies the DML property ($f$ will be an endomorphism automatically under such condition).
\end{proposition}

\begin{prf}
Just the same as the proof of Theorem 1.0.2.
\end{prf}

\begin{remark}
In fact, the endomorphism $f$ will satisfy the DML property if $\sigma f^{n}\sigma^{-1}$ satisfies for some positive integer $n$ and some automorphism $\sigma$ of $\mathbb{P}_{K}^{N}$. In particular, we can see that $f:\mathbb{P}_{K}^{N}\rightarrow\mathbb{P}_{K}^{N}$,
$$\begin{bmatrix} x_{0} \\ x_{1} \\ x_{2} \\ \vdots \\ x_{N} \end{bmatrix} \mapsto \begin{bmatrix} \sum\limits_{i=0}^{N}a_{0i}x_{i}^{q}+g_{0}(x_{0}^{p},x_{1}^{p},x_{2}^{p},\cdots,x_{N}^{p}) \\ \sum\limits_{i=0}^{N}a_{1i}x_{i}^{q}+g_{1}(x_{0},x_{1}^{p},x_{2}^{p},\cdots,x_{N}^{p}) \\ \sum\limits_{i=0}^{N}a_{2i}x_{i}^{q}+g_{2}(x_{0},x_{1},x_{2}^{p},\cdots,x_{N}^{p}) \\ \vdots \\ \sum\limits_{i=0}^{N}a_{Ni}x_{i}^{q}+g_{N}(x_{0},\cdots,x_{N-1},x_{N}^{p}) \end{bmatrix}$$
satisfies the DML property in which $K,A=(a_{ij})_{(N+1)\times(N+1)},q$ and the coefficients of $g_{0},\cdots,g_{N}$ are same as in Theorem 1.0.2 or Proposition 5.2.6 because $f^{N}$ has the expected form.
\end{remark}

\bibliographystyle{alpha}
\bibliography{reference}

\end{spacing}
\end{document}